%%%%%%%%%%%%%%%%%%%%   Geometry and Topology: 1997-1.tex  %%%%%%%%%%%%%
%%%%        
%%%%              Groups acting on Cat(0) cube complexes
%%%%             
%%%%                 Graham Niblo and Lawrence Reeves
%%%%  
%%%%                 Published in Volume 1(1997) 1--7
%%%%
%%%%                    Publication date  7-2-97
%%%%
%%%%
%%%%%%%%%%%%%%%%%%          gtplain.tem       %%%%%%%%%%%%%%%%%%
%
%  Template for articles in Geometry and Topology using plain TeX
%  and the gt macro package(s).  For instructions on using these
%  packages see  gtmacins.tex  or  gtmacins.dvi  in  gt/macros
%
%
%                  %  gtmacs  is the standard macro package.  
%                  %  If you want to use AMS fonts then comment 
%\input gtmacs     %  out this line and uncomment the following
%                  %  line (and the next if using AMS names):
\input gtmacros
\input amsnames
%
%  add further \inputs (eg  pictex  epsf  rlepsf  labelfig ) here
%
%  Original spacing (changed in gtmacros)
%
\medskipamount=6pt plus 2pt minus 2pt
\mathsurround=1pt
%
%%%%%%%%%%%%%%%%%%%%%%%%%%%%%%%%%%%%%%%%%%%%%%%%%%%%%%%%%%%%%%%%%%%%
% Preamble, our macros
%
\def\R{{\Bbb R}}                     % The real numbers
\def\QH#1{\hbox{$\Bbb {OH}^#1$}}     % Abbreviation for quaternionic hyp. space
\def\em{\bf}                         % Emphasis in definitions etc.
\def\lk#1{lk(#1)}                    % The link of a cell #1

%%%%%%%%%%%%%%%%%%%%%%%%%%%%%%%%%%%%%%%%%%%%%%%%%%%%%%%%%%%%%%%%%%%%%

%%%%%%%%%%%%          Main text starts here        %%%%%%%%%
%
%                       References first
%
\reflist
%  Type references below using the format of the following examples.
%  (Delete these examples unless you happen to need them !)
%  Don't forget to leave a blank line between references.

\refkey\B 
{\bf M Bridson},
{\it Geodesics and curvature in metric simplicial complexes},
from:\break ``Group Theory from a Geometrical Viewpoint'',  
E Ghys et al (eds.), 
World Scientific 
(1991) 
373--463

\refkey\BH 
{\bf M Bridson}, {\bf A Haefliger},
{\it Metric spaces of non-positive curvature},
in preparation

\refkey\BJS 
{\bf M Bozejko}, {\bf T Janusckiewicz}, {\bf R\,T Spatzier},
{\it Infinite Coxeter groups do not have Kazhdan's property},
J. Operator Theory 
19
(1988)
63--37

\refkey\D 
{\bf M Davis},
{\it Buildings are {\rm CAT($0$)}},
from: ``Geometric methods in group theory'',
P\,H Kropholler, G\,A Niblo and R Stohr (eds.),
LMS Lecture Note Series,
Cambridge University Press

\refkey\G 
{\bf M Gromov},
{\it Hyperbolic groups},
from: ``Essays in group theory'',
S\,M Gersten (ed.),
MSRI Publ. 8,
Springer--Verlag
(1987)
75--267

\refkey\HV 
{\bf P de la Harpe}, {\bf A Valette},
{\it La propri\'et\'e {\rm (T)} de Kazhdan pour les groupes localement compacts},
Asterisque 175 (1989),
Soci\'et\'e Math\'ematique de France

\refkey\NR 
{\bf G\,A Niblo}, {\bf L\,D Reeves},
{\it Coxeter groups act on {\rm CAT($0$)} cube complexes},
preprint

\refkey\S 
{\bf M Sageev},
{\it Ends of group pairs and non-positively curved cube complexes},
Proc. London Maths. Soc. (3) 71 (1995)
585--617

\endreflist
%
%                          Title page  
\input gtoutput

\volumenumber{1}\papernumber{1}
\volumeyear{1997}\pagenumbers{1}{7}
\published{7 February 1997}

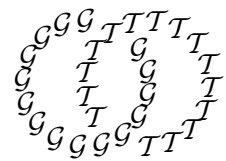
\title{Groups acting on CAT($0$) cube complexes}           
\shorttitle{Groups acting on CAT(0) cube complexes}           
\authors{Graham Niblo\\Lawrence Reeves}

\address{Faculty of Mathematical Studies\\University of Southampton\\%
Highfield\\Southampton, SO17 1BJ, UK}
\secondaddress{Institute of Mathematics\\Hebrew University of 
Jerusalem\\Givat Ram\\Jerusalem 91904, Israel}
\email{gan@maths.soton.ac.uk\\lawrence@math.huji.ac.il}

\abstract

We show that groups satisfying Kazhdan's
property (T) have no unbounded actions on finite dimensional CAT($0$) cube
complexes, and deduce that there is a locally CAT($-1$) Riemannian
manifold which is not homotopy equivalent to any finite dimensional, 
locally CAT($0$) cube
complex.

\endabstract




\proposed{Walter Neumann}\received{28 October 1996}
\seconded{David Gabai, Robion Kirby}\accepted{6 February 1997}
\maketitlepage

%%%%%%%%%%%%%%%%%%%%   End of title page
%
%%%%%%%%%%%%%%%%%%%    Start of main body of article

\section{Introduction}
%%%%%%%%%%%%%%%%%%%%%%%%%%%%%%%%%%%%%%%%%%%%%%%%%%%%%%%%%
%%%%%% Intro to CAT($0$) spaces
%%%%%%%%%%%%%%%%%%%%%%%%%%%%%%%%%%%%%%%%%%%%%%%%%%%%%%%%%
The CAT($\chi$) inequality gives a measure of the curvature of a geodesic metric
space
$X$ by comparing the width of the geodesic triangles in $X$ with those of the
corresponding triangles in the simply connected Riemannian manifold of constant
curvature $\chi$. The theory of CAT($\chi$) metric spaces is described in
[\BH].

A geodesic metric space
$X$ is said to be locally CAT($\chi$) if around every point there is a non-empty
ball which is geodesically convex and 
CAT($\chi$). It follows from the definition that
$X$ is locally CAT($0$) if every point has a CAT($\chi$)
neighbourhood for some $\chi$ less than or equal to $0$, where $\chi$ is
allowed to vary with the point. 

Examples of locally CAT($0$) spaces are furnished by Riemannian manifolds
whose sectional curvatures are all bounded above by
$0$, by what Davis calls in [\D] the ``(correctly defined) geometric
realization" of
a Tits' building, and by compact piecewise Euclidean polyhedral complexes
satisfying certain curvature conditions on the link of every cell, as
studied by Bridson in his thesis, [\B]. If the cells of such a complex are all
isometric to Euclidean cubes then the curvature condition can be expressed in
terms of the combinatorics of the links of the cells (Gromov's ``no bigon''
and ``no triangle'' conditions [\G], outlined here in  section 3). It
is an open problem to decide which CAT($0$) spaces can be given a piecewise
Euclidean polyhedral structure which is also CAT($0$).

%%%%%%%%%%%%%%%%%%%%%%%%%%%%%%%%%%%%%%%%%%%%%%%%%%%%%%%%%
%%%%%% Intro to the main theorem
%%%%%%%%%%%%%%%%%%%%%%%%%%%%%%%%%%%%%%%%%%%%%%%%%%%%%%%%%
In this paper we establish the following theorem:

\proclaim{Theorem A} For any $n\geq 2$ there is a smooth locally {\rm
CAT($-1$)} Riemannian manifold
$M$ of dimension $4n$ which is not homotopy equivalent to any finite
dimensional, locally {\rm CAT($0$)} cube complex.
\endproc\ppar

The manifold $M$ is a quotient of quaternionic hyperbolic
space \QH n (of real dimension $4n$) by a
discrete group $G$ of isometries. The isometry group of \QH n
is the symplectic group $Sp(n,1)$; in fact we may choose
$G$ to be any discrete co-compact lattice in $Sp(n,1)$, and let $M$
be the corresponding quotient of \QH n, so $G=\pi_1(M)$.
Now if
$M$ was homotopy equivalent to a
compact, locally CAT($0$) cube complex
$X$, then $G$ would act freely and properly discontinuously on the
universal cover $\tilde{X}$ which
is a non-compact finite dimensional CAT($0$) cube complex (see [\B] for a
discussion of how the local curvature properties of $X$ give global curvature
constraints on the universal cover). We will see that
$G$ admits no such actions; indeed we will show
that any cellular action of $G$ on a finite dimensional CAT($0$) cube
complex must have  a global fixed point. This will follow from the fact
that $G$ has Kazhdan's property
(T) [\HV], together with the following result:
%%%%%%%%%%%%%%%%%%%%%%%%%%%%%%%%%%%%%%%%%%%%%%%%%%%%%%%%%
%%%%%% Intro to Kazhdan's property (T)
%%%%%%%%%%%%%%%%%%%%%%%%%%%%%%%%%%%%%%%%%%%%%%%%%%%%%%%%%
\proclaim{Theorem B} If $G$ is a group satisfying Kazhdan's property 
{\rm (T)} and
$X$ is a finite dimensional {\rm CAT($0$)} cube complex on which $G$ 
acts cellularly
(and therefore isometrically), then the action has a global fixed
point.
\endproc\ppar

We note that since $M$ is a compact locally CAT($-1$) space, $\pi_1(M)$
is a
$\delta$--hyperbolic group, and so this gives an example of a
$\delta$--hyperbolic group which cannot act without fixed points on
any finite dimensional CAT($0$) cube complex.

Further examples of
$T$--groups are furnished by  groups acting simply transitively on Tits' buildings 
modelled on the $\tilde{A_2}$ simplex, ie the Euclidean triangle with angles   $\pi/4$, $\pi/4$ and $\pi/2$, [\HV]. As Davis shows
in [\D], these buildings are CAT($0$), so their quotients are locally CAT($0$). Their quotients give  examples of  locally CAT($0$)
simplicial complexes which are not homotopy equivalent to any finite dimensional, locally CAT($0$) cube complexes.

%%%%%%%%%%%%%%%%%%%%%%%%%%%%%%%%%%%%%%%%%%%%%%%%%%%%%%%%%
%%%%%% Intro to Coxeter groups
%%%%%%%%%%%%%%%%%%%%%%%%%%%%%%%%%%%%%%%%%%%%%%%%%%%%%%%%%
The proof of Theorem B is modelled on (or stolen from, depending
on your point of view) the proof by   Bozejko et al
[\BJS] that finitely generated Coxeter groups do not have
property
(T).
%%%%%%%%%
%%%%%Property $T$ groups are fg so true for all Coxeter groups
%%%%%%%%%
It can be shown [\NR] that finitely generated Coxeter groups
act effectively, properly discontinuously and cellularly on  finite
dimensional CAT($0$) cube complexes so Theorem A may be viewed as a
generalisation of their result.

%The relationship between the cubing and the affine construction of
%the Coxeter complex? Seeing the null vectors as the extra dimensions?

%%%%%%%%%%%%%%%%%%%%%%%%%%%%%%%%%%%%%%%%%%%%%%%%%%%%%%%%%
%%%%%% The layout of the paper
%%%%%%%%%%%%%%%%%%%%%%%%%%%%%%%%%%%%%%%%%%%%%%%%%%%%%%%%%
In section 2 we will recall the definition of a  ``conditionally
negative kernel'' on a group as given in [\HV] and the
characterisation of Kazhdan's property
(T) in terms  of these kernels. We will then outline the proof of Theorem
B. In section
3 we will give the definition of a CAT($0$) cube complex and derive the
technical result we need to complete the missing steps of the proof.

\section{Kazhdan's property (T)}

%%%%%%%%%%%%%%%%%%%%%%%%%%%%%%%%%%%%%%%%%%%%%%%%%%%%%%%%%
%%%%%% Definition of property (T) and outline proof of main theorem
%%%%%%%%%%%%%%%%%%%%%%%%%%%%%%%%%%%%%%%%%%%%%%%%%%%%%%%%%
In [\HV] de la Harpe and Valette give several equivalent
formulations of Kazhdan's property (T) for locally compact groups.
Although it is not the most intuitive one we will use the following
criterion to establish our main result:

\rk{Definition}A {\em conditionally negative kernel} on a set
$V$ is a function $f\co V\times V\longrightarrow\R$ such that for any
finite subset $\{v_1,\ldots,v_n\}\subset V$ and any real numbers
$\{\lambda_1,\ldots,\lambda_n\}$ such that
$\sum\limits_{i}\lambda_i=0$ the following inequality holds:
$$
\sum\limits_{i,j}\lambda_i\lambda_j f(v_i,v_j)\leq 0
$$
A  conditionally negative kernel on a group
$G$ is a conditionally negative kernel on the set of elements of $G$
such that for any $g,h,k$ in $G$,
$f(gh,gk)=f(h,k)$.

According to [\HV] a finitely
generated  group
$G$ has Kazhdan's property (T) (or  is  a {\em $T$--group}) if and only if every
conditionally  negative kernel on
$G$ is bounded.

Our goal is to prove the following theorem. The statement is followed by an
outline
proof and the missing lemma is given in the following section.

\proclaim{Theorem B} If $G$ is a group satisfying Kazhdan's property {\rm (T)}
and
$X$ is a finite dimensional {\rm CAT($0$)} cube complex on which $G$ acts cellularly
(and therefore isometrically), then the action has a global fixed
point.
\endproc

\prf In section 3 we will
show  that for any CAT($0$) cube complex $X$, the simplicial  metric $D$
on the 1--skeleton of $X$ restricts to a conditionally negative kernel on the
vertex set $X^{(0)}$. Since $G$ acts cellularly on $X$, $D$ is invariant under
the action, so setting $f(g,h)=D(gv,hv)$ gives a
conditionally
negative kernel on the group $G$. Since
$G$ has Kazhdan's property (T) the conditionally negative kernel must be
bounded so the
orbit $Gv$ is bounded in the metric $D$ and therefore, since $X$ is finite
dimensional,  in the metric $d$ (see section 3).
Since
$X$ is finite dimensional it is complete [\B] and in a complete CAT($0$) metric
space any isometric action with a bounded orbit has a global fixed point
([\B] again),
completing the proof of Theorem B.
\qed

%%%%%%%%%%%%%%%%%%%%%%%%%%%%%%%%%%%%%%%%%%%%%%%%%%%%%%%%%%%%%%%%%%%%%%%%%%%%
%%%%%%%%%%%%%%%%%%%%%%%%%%%%%%%%%%%%%%%%%%%%%%%%%%%%%%%%%%%%%%%%%%%%%%%%%%%%
\section{CAT($0$) cube complexes}

A {\em cube complex} $X$ is a
metric polyhedral complex in which each cell is isometric to the Euclidean
cube $[-1/2,1/2]^n$, and the gluing maps are isometries. If there is a bound on
the dimension of the cubes then the complex carries a complete geodesic
metric, [\B].

\rk{Definition} A cube complex is {\em non-positively curved} if for
any cube $C$ the following conditions on the link of $C$, $\lk{C}$, are
satisfied:
\items
\item{({\it i\/})}({\em no bigons})\stdspace For each pair of vertices in
$\lk{C}$ there is at most one edge containing them.
\item{({\it ii\/})}({\em no triangles})\stdspace Every
edge cycle of length three in $\lk{C}$ is contained in a 2--simplex of $\lk{C}$.
\enditems

The
following theorem of Gromov relates the combinatorics and the geometry of the
complex.

\proclaim{Lemma} {\rm (Gromov, [\G])}\stdspace
A cube complex
$X$ is locally {\rm CAT($0$)} if and only if it is non-positively curved, and it
is {\rm CAT($0$)} if and only if it is non-positively curved and simply connected.
\endproc

\rk{Examples}\rm  Any graph may be regarded as a 1--dimensional cube
complex, and the curvature conditions on the links are trivially satisfied.
The graph is CAT($0$) if and only if it is a tree. Euclidean space also has
the structure of a CAT($0$) cube complex with its vertices at the integer
points.
\ppar

A {\em midplane} of a cube $[-1/2,1/2]^n$ is its intersection with a
codimension 1 coordinate hyperplane. So every $n$--cube contains $n$
midplanes each of which is an $(n-1)$--cube, and any $m$ of which intersect in
a $(n-m)$--cube. Given an edge in a non-positively curved cube complex, there
is a unique codimension 1 hyperplane in the complex which cuts the edge
transversely in its midpoint. This is obtained by developing the coordinate
hyperplanes in the cubes containing the edge. In the case of a tree the
hyperplane is the midpoint of the edge, and in the case of Euclidean space it
is a geometric hyperplane.

In  general the hyperplane is analogous to an immersed codimension 1
submanifold in a  Riemannian
manifold and the immersion is actually a local isometry. An application of
the Cartan--Hademard
theorem, [\G, section 4], then shows that the hyperplane is isometrically
embedded, and
furthermore any hyperplane in a CAT($0$) cube complex separates it into two
components referred to as the half spaces associated with the hyperplane.
This is a consequence of the fact that the complex is simply connected. The
hyperplane gives rise to a 1--cocycle which is necessarily trivial, and hence
the hyperplane separates the space.

The set of vertices of the cubing $X$ can be viewed as a discrete
metric space, where the metric $D(u,v)$ is given by the length of
the shortest edge path between the
vertices $u$ and $v$ in the 1--skeleton of $X$. If $X$ is finite dimensional
this metric is quasi-isometric to the CAT($0$) metric (it is at most
$\sqrt{n}$ times the CAT($0$) metric, where $n$ is the dimension of the
complex; this does not require the complex to be cocompact).
Sageev, [\S], observed that the shortest path in the 1--skeleton crosses any
hyperplane at most once, and since every edge crosses exactly one
hyperplane,
the distance between two vertices is the number of hyperplanes separating
them. Recalling that a hyperplane separates the vertices into two half
spaces $U^+$ and $U^-$, it follows that
$$D(u,v)=\sum\limits_U \chi_U(u)(1-\chi_U(v))$$
where $\chi_U$ is the characteristic function of the half space $U$, and $U$
ranges over all the half spaces in $X$. Although there are infinitely many half
spaces, the fact that only finitely many of them separate any pair of vertices
means that this will be a finite sum.

Now suppose that $G$ is a group acting cellularly on $X.$
The action restricts to the 1--skeleton and therefore preserves the metric
$D,$ so we have established parts 2 and 3 of the following lemma:

\proclaim{Technical Lemma} For any {\rm CAT($0$)} cube
complex $X$ the simplicial metric $D$ on the 1--skeleton of $X$ satisfies the
following properties:
\items
\item{1)}$D$ is a conditionally negative kernel on the vertex set of $X$.
\item{2)}$D$ is invariant under the action of $G$ on the vertex set.
\item{3)}For any vertices
$u, v\in X$,
$d(u,v)\leq D(u,v)\leq
\sqrt{n}d(u,v)$ where $n$ is the dimension of the cube complex.
\enditems
\endproc

\prf
It remains to establish that the metric $D$ is a conditionally negative kernel.
Let
$u_1,u_2,\ldots, u_n\in G$ and
$\lambda_i\in\R$ with
$\sum\limits_{i}\lambda_{i}=0$. Then
$$
\sum\limits_i\sum\limits_j\lambda_i\lambda_jD(u_i,u_j)=\sum\limits_i\sum\limits_
j
\lambda_i\lambda_j\sum\limits_U\chi_U(u_i)
\left(1-\chi_U(u_j)\right)
$$
where the third sum can taken over the finitely many half
spaces
$U$ separating the vertices $\{u_1,\ldots, u_n\}$ (the other half spaces
all contribute 0 to the sum). This triple sum can now be rearranged to get
$$
\sum\limits_i\sum\limits_j\lambda_i\lambda_j\sum\limits_U\chi_U(u_i)-
\sum\limits_i\sum\limits_j\lambda_i\lambda_j\sum\limits_U\chi_U(u_i)\chi_U(
u_j)
$$
The first term is zero since it can be rearranged to give
$\sum\limits_j\lambda_j\sum\limits_i\lambda_i\sum\limits_U\chi_U(u_i)$, and
\hbox{$\sum\limits_j\lambda_j=0$} by hypothesis. Likewise the second term
can be
rewritten as
$$
\sum\limits_U\sum\limits_i\lambda_i\chi_U(u_i)\sum\limits_j\lambda_j\chi_U
(u_j)
=\sum\limits_U\left(\sum\limits_i\lambda_i\chi_U(u_i)\right)^2
$$
which is positive. It follows that the entire sum is negative as required.
\qed

\np

%%%%%%%%%%%%%%%%%%%%  End of main body of article
%
\references         %  Print the references here
%
%    Any appendix material (to go after the references) here 
%
\bye